\documentclass[12pt]{article}

\usepackage{leqno}

\usepackage{amssymb}
\usepackage{amsmath}

 \newcommand{\R}{{\mathbb R}}
 
 \newcommand{\CQFD}
 {
 \mbox{}
 \nolinebreak
 \hfill
 \rule{2mm}{2mm}
 \medbreak
 \par
 }

 \newtheorem{Theo}{Theorem}

 \newtheorem{lemma}{Lemma}

 \newtheorem{Coro}{Corollary}
 \newtheorem{Remark}{Remark}

\begin{document}



\title{
{\Large$L_1^{}$}-optimal linear programming estimator
for periodic frontier functions with H\"older continuous derivative
}

\author{A.V.~Nazin.\thanks{The work of A.V.~Nazin was carried out during his stay
in MISTIS Project, Inria Grenoble Rh\^{o}ne-Alpes, June and October
2013;
partially supported by PreMoLab/MIPT,
RF government grant \mbox{№\,11.G34.31.0073}.
}\\
(Trapeznikov Institute of Control Sciences RAS, Moscow, Russia),\\
S.~Girard.
\\
(LJK, Inria Grenoble Rh\^one-Alpes, 
Grenoble, France%
)%
}
\date{}
\maketitle

\begin{abstract}
We propose a new estimator based on a linear programming method
for smooth frontiers of sample points.
The derivative of the frontier function is
supposed to be H\"older continuous.
The estimator is defined as a linear combination of kernel functions being sufficiently regular,
covering all the points and whose associated support is of smallest surface.
The coefficients of the linear combination are
computed by solving a linear programming problem.
The $L_1$ error between the estimated and the true frontier functions
is shown to be almost surely converging to zero, and the rate of convergence
is proved to be optimal.
\end{abstract}



\section{Introduction}
\label{secIntro}

Many proposals are given in the literature for estimating a set $S$
given a finite random set of points drawn from the interior.
Here, we focus on the case where the unknown support can be written
as $S=\{(x,y):\;0\leq x\leq 1~;~0\leq y \leq f(x)\}$,
where $f$ is an unknown function.
The initial problem reduces to estimating $f$, called the frontier or the boundary,
from random pairs $(X,Y)$ included in $S$.

Under monotonicity assumptions, the frontier can also be interpreted
as the endpoint of $Y$ given $X\leq x$. Specific estimation techniques have
been developed in this context, see for instance {\sc Deprins} {\it et al.}~\cite{DepSimTul},
{\sc Farrel}~\cite{Farrel}, {\sc Gijbels} {\it et al.}~\cite{Gijbels2}. We also refer to
{\sc Aragon} {\it et al.}~\cite{Aragon}, Cazals {\it et al.}~\cite{Cazals}, {\sc Daouia} \& {\sc Simar}~\cite{Daouia} for the definition of robust estimators.

In the general case, that is without monotonicity assumptions, {\sc
Girard} \& {\sc Jacob}~\cite{GirJac4} introduced an estimator based
upon kernel regression on high power-transformed data. In the
particular case where $Y$ given $X=x$ is uniformly distributed they
proved that this estimator is asymptotically Gaussian with the
minimax rate of convergence for Lipschitzian frontiers.
\textrm{(Loosely speaking, under the rate of convergence we understand
infinitely small positive number sequence which
characterizes the convergence to zero of a norm of the estimation error,
as the sample size $ N \to\infty $.)}
 Compared to the extreme-value based estimators ({\sc Geffroy}~\cite{Geff1},
{\sc Girard} \& {\sc Jacob}~\cite{GirJac,GirJac2,GirJac3}, {\sc
Girard} \& {\sc Menneteau}~\cite{GirMen}, {\sc H\"ardle} {\it et
al.}~\cite{Hardle2}, {\sc Menneteau}~\cite{Mennet}), projection
estimators ({\sc Jacob} \& {\sc Suquet}~\cite{JacSuq}), or piecewise
polynomial estimators ({\sc Hall} {\it et al.}~\cite{Hall2}, {\sc
Knight}~\cite{KK}, {\sc Korostelev} \& {\sc
Tsybakov}~\cite{KorTsy2}, {\sc Korostelev} {\it et
al.}~\cite{KorTsy}, {\sc H\"ardle} {\it et al.}~\cite{Hardle}), this
estimator does not require a partition of the support $S$. When the
conditional distribution of $Y$ given $X$ is not uniform, this
estimator is still convergent ({\sc Girard} \& {\sc
Jacob}~\cite{GirJac4}, Theorem~1) but may suffer from a strong bias
({\sc Girard} \& {\sc Jacob}~\cite{GirJac4}, Table~1). A
modification of this estimator has been proposed by {\sc Girard}
{\it et al.}~\cite{Stu1,Stu2} to tackle the situation where the
conditional distribution function of $Y$ given $X=x$ decreases at a
polynomial rate to zero in the neighborhood of the frontier $f(x)$.
The asymptotic normality as well as the strong consistency of the
estimator are established.

The estimator proposed in {\sc Bouchard} {\it et al}~\cite{BGIN}
for estimating $S$ shares some common characteristics with the one of
{\sc Girard} \& {\sc Jacob}~\cite{GirJac4}. It assumes that $Y$ given $X=x$ is uniformly distributed
but does not require a partition of the support. Besides,
it is defined as a kernel estimator obtained by smoothing
some selected points of the sample.
These points are, however,
chosen automatically by solving a linear programming problem
to obtain an estimated
support covering all the points and with smallest surface.
From the theoretical point of view, this estimator is shown
to be consistent for the $L_1$ norm. An improvement of this
estimator has been proposed in {\sc Girard} {\it et al.}~\cite{GIN}
in order to reach the optimal minimax $L_1$ rate of convergence (up to a logarithmic factor)
for Lipschitzian frontiers.

In this paper, we propose an adaptation of these methods for estimating
smoother frontiers: It is assumed that the first derivative of frontier
is \textrm{H\"older continuous}.
The resulting estimator is proved
to reach the optimal minimax $L_1$ rate of convergence (up to a logarithmic factor).
The paper is organized as follows. The estimator is defined in Section~\ref{subsecLPprobl}.
Assumptions and preliminary results are given in Section~\ref{preliresult}
while our main result is established in Section~\ref{secmain}. Proofs are postponed to the
Appendix.

\section{Problem statement and boundary estimator}
\label{subsecLPprobl}

Let all the random variables be defined on a probability space
$(\Omega,\mathcal{F},P)$. The problem under consideration is to
estimate an unknown
$1$-periodic function $f: \mathbb{R}\to
(0,\infty)$, that is $f(x+1)=f(x)$ for all $x\in\R$, on the basis of
 independent observations $(X_i,Y_i)_{i=\overline{1,N}}$
uniformly distributed in
  \begin{equation}\label{defS1}
  S \triangleq \{(x,y) : \, 0  \leq x\leq 1\,,\,\, 0\leq y\leq f(x) \}\,.
 \end{equation}
Note that, the notation $i=\overline{m,n}$ is used for $i=m,\dots,n$.
Since $f$ is $1$-periodic, it is convenient to
extend the indices of data $(X_i,Y_i)$ out of those of
$i=\overline{1,N}$ by periodic continuation w.r.t. $x$. Therefore,
we put $(X_i,Y_i)=(X_{i+N}-1,Y_{i+N})$ for $i=\overline{1-N,0}$, and
$(X_i,Y_i)=(X_{i-N}+1,Y_{i-N})$ for $i=\overline{N+1,2N}$.
\begin{Remark}\label{rem:1}
\textrm{An example of the boundary of $ 2 \pi $-periodic function occurs in the description
of the boundary of a convex planar set in polar coordinates with the center inside the set.
The normalization of the polar angle allows to come to $ 1 $-periodic function.}
\end{Remark}
\begin{Remark}\label{rem:2}
\textrm{Note that the condition of periodicity of the boundary function is, on the one hand,
a significant assumption that distinguish a special class of problems,
on the other hand, the technical condition that allows (together with the conditions on the kernel function, see below)
to avoid the ``difficulties at the borders'' of the interval $[0,1]$,
simplify the calculations in the analysis of the estimation error and arrive to the optimal rate of convergence.
}
\end{Remark}
Letting
\begin{equation}\label{eq:Cf}
    C_f \triangleq \displaystyle\int_0^1 f(u)\,du 
    \,,
\end{equation}
each variable $X_i$ is distributed in $[0,1]$ with p.d.f.
$f(\cdot)/C_f$ while $Y_i$ has an uniform conditional distribution
with respect to $X_i$ in the interval $[0,f(X_i)]$.
In what follows, it is assumed that $f\in\Sigma(1,\beta,L^{}_{\beta}\,)$,
$1<\beta\leq 2$, i.e. the function $f:{\mathbb R}\to(0,\infty)$ is
$1$-periodic, continuously differentiable
with
H\"older continuous derivative $f^\prime$ having exponent
$\beta-1$
and upper bound for H\"older coefficient $L^{}_{\beta}\,$\,:
\begin{equation}\label{LipschitzF}
|f'(x)-f'(u)|\leq L^{}_{\beta}\,|x-u|^{\beta-1}\quad\forall\,
x,u\in {\mathbb R}\,.
\end{equation}
\noindent The considered estimator
$\widehat{f}_N:[0,1]\to[0,\infty)$ of the frontier is chosen from
the family of functions \noindent
\begin{equation}\label{estimator}
\begin{array}{l}
\widehat{f}_N(x) = \displaystyle \sum_{i=1}^N \alpha_i\,
K_h(x,X_i),\qquad \alpha_i \geq 0,\quad i=1,\dots,N,
\end{array}
\end{equation}
with kernel function
\begin{equation}\label{estimatorsuite}
K_h(x,t)= \frac{1}{h}\,K\left(\frac{x-t}{h}\right) + \left|
\begin{array}{ll}
 0 &  \mbox{ if }\, h<t<1-h,\\
\displaystyle
 \frac{1}{h}\;K\!\left(\frac{x-t-1}{h}\right) & \mbox{ if }\, 0\leq t\leq h,
 \smallskip\\
\displaystyle
 \frac{1}{h}\;K\!\left(\frac{x-t+1}{h}\right) & \mbox{ if }\, 1-h\leq t\leq 1,
\end{array}
\right.
\end{equation}
being defined for all $(x,t)\in[0,1]^2$
 and
$h\in(0,1/2)$,
where $K$ is a given sufficiently smooth centered density function\,
$K:\R\to[0,\infty)$ with support included in $[-1,1]$, see assumption B2, Section~\ref{preliresult}.
The bandwidth parameter $h$
depends on $N$ such that $h\to 0$ as $N\to\infty$.
\medskip

\begin{Remark}\label{rem:nonperiodic}
The estimate (\ref{estimator}) is not supposed to be periodic
itself.
\end{Remark}
\begin{Remark}\label{Remark1}
Since the supports of the different terms appearing in
(\ref{estimatorsuite}) do not intersect, the kernel may
be rewritten as
\begin{equation}\label{estimatorsuitePlus3}
K_h(x,t)= \sum_{j=-1}^{1}
            \frac{1}{h}\,K\left(\frac{x-t+j}{h}\right)
%
\end{equation}
for $(x,t)\in[0,1]^2$.
The kernel function~(\ref{estimatorsuitePlus3}) is thus as smooth w.r.t. $t$ as the density
function $K(\cdot)$ is. For instance, one always has
$K_h(x,t)\leq K_{\max}/h$ and the $k$-th derivative bound
$|\partial^{\,k\!}K_h(x,t)/\partial t^k| \leq
K_{\max}^{(k)}/h^{k+1}$ if $K(\cdot)$  has a continuous $k$-th
derivative.
\end{Remark}
\begin{Remark}\label{rem:5}
\textrm{The optimal choice of the bandwidth $h$ is carried out here in the spirit of \cite{GIN}.
Also see Remark~\ref{rem:OptimLowerBound} below at the end of Section~\ref{subsec:LB}.}
\end{Remark}
\begin{Remark}\label{rem:WiderInterval}
In the definition of estimator (\ref{estimator}),
the kernel~(\ref{estimatorsuite}) is introduced for  $x\in[0,1]$.
However, it is convenient to introduce a wider interval like
$[-h,1+h]$ for variables $x$ and $t$ in the kernel function
(\ref{estimatorsuite}) and define additional points $X_{-}\in[-h,0]$
and $X_{+}\in[1,1+h]$ a.s., see below (\ref{eq:X0XN1}).
\end{Remark}

\noindent As it is proved below in Lemma \ref{LemSurfApprox} the
surface of the estimated support
  \begin{equation}\label{defSesti}
  \widehat S_N \triangleq \{(x,y) : \, 0  \leq x\leq 1\,,\,\, 0\leq y\leq \widehat f_N(x) \}
 \end{equation}
 is given by
\begin{equation}\label{SurfApprox}
\int_0^1 \widehat{f}_N(x)\, dx = \sum_{i=1}^N \alpha_i \,.
\end{equation}
This suggests to define the estimator of the parameter vector
$\alpha=(\alpha_1,\dots,\alpha_N)^T$ as a solution of the following
optimization problem
\begin{eqnarray}\label{IPgoalS}
J_P^* &\triangleq& \min_\alpha \sum_{i=1}^N \alpha_i
\\
\mathrm{subject~to,} && \mathrm{for~all~}i=\overline{1,N}, \nonumber
\\
&& \widehat{f}_N(X_i)
                    +\left(X_j - X_i\right)\widehat{f}^{\;\prime}_N(X_i)
 \geq Y_j\,,\quad\forall\,j:\,|X_j - X_i|\leq h\,,
\label{constr1S}\\
&& |\widehat{f}^{\;\prime\prime}_N(X_i)|
 \leq
2L^{}_{\beta}\, K''_{\max}\, \frac{\log N}{Nh^3}\,,
 \label{constrPrime2S}\\
&&\sum_{i=1}^N \alpha_i \, \mathbf{1}\{ (m-1)/m_h\leq X_i < m/m_h\}
\leq C_\alpha h \,,\quad m=\overline{1,m_h} \,,
 \label{constrSumAlp} \\
&&0\leq\alpha_i\,,
\label{constr2S}
\end{eqnarray}
where $m_h=\lfloor{1/h}\rfloor$ is the integer part of $1/h$,
$$
K^{''}(x,u) \triangleq \frac{\partial^{\,2}}{\partial
x^2}\,K_{h}(x,u)\,, \quad (x,u)\in [0,1]^2\,,
$$
and $\mathbf{1}\{\cdot\}$ is the
indicator function 
which
equals $1$ if the argument condition holds true, and $0$ otherwise.
The value of the positive parameter $C_\alpha$ in the constraints~(\ref{constrSumAlp})
is discussed in Section~\ref{secmain}.
Evidently, this optimization problem
represents a linear program (LP).
Therefore, we call the defined boundary estimator as the LP-estimator (\ref{IPgoalS})--(\ref{constr2S}).
\section{Basic assumptions and preliminary results}\label{preliresult}

The basic assumptions on the unknown boundary function
$f:\R\to(0,+\infty)$ are:
\begin{itemize}
\item[A1.] $f(x)=f(x+1)$ and
$0< f_{\min} \leq f(x) \leq f_{\max} <\infty$, for all $x\in \R$.
 \item[A2.] $f(x)$ is
 continuously differentiable
having the
H\"older exponent $\beta-1$ for function derivative $f^\prime$,
i.e.
$$|f'(x)-f'(y)|\leq L^{}_{\beta}\, |x-y|^{\beta-1}
 \mathrm{~~~for~all~~~} x,y \in [0,1]\,,
$$
where constants $\, L^{}_{\beta}\,<\infty$ and $\beta\in(1,2]$ are
supposed to be given.
\end{itemize}
\medskip

\noindent The following assumptions on the kernel function are
introduced:
\begin{itemize}
  \item[B1.] $K:\R\to[0,+\infty)$ has a compact support:
   $\displaystyle{\mathop{\mathrm{supp\,}}_{t\in\R}}K(t)=[-1,1]$,
  \item[B2.] $\displaystyle \int_{-1}^1 K(t)\, dt =1 $,\,\,
   $\displaystyle \int_{-1}^1 t\, K(t)\, dt =0 $,
  \item[B3.] $K$ is four times continuously differentiable.
\end{itemize}
The next two lemmas are of analytical nature. They can be interpreted
as extensions of Bochner's Lemma for controlling the smoothing error
introduced by the kernel.
%

\begin{lemma}\label{lem:IntegrKer}
Let $f$ be a $1$-periodic function on $\mathbb{R}$. Then, assumptions
B1 and B2 imply
\begin{equation}\label{eq:LemKerApprox1}
    \int_0^1 f(u) K_h(x,u)\,du = \int_{-1}^1 f(x-hv) K(v)\,dv
     \qquad \forall\,x\in[0,1]\,.
\end{equation}
Particularly, for $f(x)\equiv 1$ one obtain
\begin{equation}\label{eq:IntUKer}
     \int_0^1 K_h(x,u)\,du
     = 1
     \qquad \forall\,x\in[0,1]\,.
\end{equation}
In addition,
\begin{eqnarray}\label{eq:IntU-XdUKer}
     \int_0^1 (u-x) K_h(x,u)\,du
     &=& 0
     \qquad \forall\,x\in[h,1-h]\,,\\
\label{eq:IntXKer1}
     \int_0^1 K_h(x,u)\,dx
     &=& 1
     \qquad \forall\,u\in[0,1]\,,\\
\label{eq:IntXUdxKer}
     \int_0^1  (x-u) K_h(x,u)\,dx
     &=& 0
     \qquad \forall\,u\in[0,1]\,.
\end{eqnarray}
\end{lemma}
Now, we quote several preliminary results on the estimator $\widehat
f^{}_N$.  First, Lemma~\ref{LemSurfApprox} establishes that the surface of
the related estimated support $\widehat S^{}_N$ equals $\sum_{i=1}^N
\alpha_i$. 
Second, functions ${\widehat f}^{}_N$ and ${\widehat
f}_N^{\,\prime}$ are proved to be Lipschitzian, see Lemma~\ref{RegLemma1} and Lemma~\ref{RegLemma} respectively.
Proofs are postponed to
Subsection~\ref{subsec:proofpreli}.


\begin{lemma}\label{LemSurfApprox}
Suppose B1, B2 are verified and $0<h<1/2$. Then, the surface of the
estimated support (\ref{defSesti}) is
\begin{equation}\label{SurfApproxUB}
 \int_0^1 \widehat{f}_N(x)\, dx =  \sum_{i=1}^N \alpha_i\,.
\end{equation}
\end{lemma}

Introduce additionally
\begin{eqnarray}\label{eq:X0XN1}
X_{-}&\triangleq& \max_{i=\overline{1,N}}\{X_i\}-1\,, \quad
X_{+}\triangleq \min_{i=\overline{1,N}}\{X_i\}+1\,.
\label{X0X1defs}
\end{eqnarray}

\begin{lemma}\label{RegLemma1}
Suppose A1, A2, and B1, B3 are verified. Then, the Lipschitz constant of the LP-estimator
(\ref{IPgoalS})--(\ref{constr2S}) is bounded by
\begin{eqnarray}\label{eq:RegLemmaP1}
L_{\widehat{f}_N} &\triangleq&
    \max_{x\in[X_{-},X_{+}]} |\widehat{f}_N^{\;\prime}(x)|
    \leq  3 C_{\alpha}\, K'_{\max}\, h^{-1}\,.
\end{eqnarray}
\end{lemma}
\begin{lemma}\label{RegLemma}
Suppose A1, A2, and B1, B3 are verified.
Moreover, let $h\to 0$ as $N\to\infty$ such that
\begin{equation}\label{RegLem8}
\lim_{N\to\infty} \frac{\log N}{N h} =0 \,.
\end{equation}
Then, there exists almost surely finite $N_4=N_4(\omega)$ such that
for any $N\geq N_4$
the 
Lipschitz
 constant for
the derivative estimator $\widehat{f}_N^{\,\prime}$ over 
interval $[X_{-},X_{+}]\supseteq[0,1]$ is bounded as follows:
\begin{eqnarray}\label{eq:RegLemmaP}
L_{\widehat{f}_N^{\,\prime}} &\triangleq&
    \max_{x\in[X_{-},X_{+}]} |\widehat{f}_N^{\;\prime\prime}(x)|
    \leq
 4\,L^{}_{\beta}\, K''_{\max}\, \frac{\log N}{Nh^3}\,.
\end{eqnarray}
\end{lemma}

\noindent 
Below in the next section (see also the proof of Theorem~\ref{Th1}),
it appears that the LP-estimator $\widehat{f}_N$
solution to the optimization problem
(\ref{IPgoalS})--(\ref{constr2S})
 defines the kernel estimator of the
support covering all the points $(X_i, Y_i)$ and having the smallest
surface. Moreover, constraints (\ref{constrPrime2S})--(\ref{constr2S})
impose  the first derivative $\widehat{f}_N^{\,\prime}$ of
the estimator to be Lipschitzian
with a particular Lipschitz constant
$L_{\widehat{f}_N^{\,\prime}}$\, given in Lemma~\ref{RegLemma}.
The constraint (\ref{constr1S}) says that, for any $X_i$,
the local linear estimate function $\widehat{f}_N(X_i) +(x-X_i)\widehat{f}_N^{\,\prime}(X_i)$
covers all points $(X_j,Y_j)$ with $x=X_j$ from interval $\{x: |x-X_i|\leq h\}$.
Additionally, the constraints $\alpha_i\geq 0$ for all
$i=\overline{1,N}$ ensure that $\widehat{f}_N(x) \geq 0$ for all
$x\in[0,1]$ since the density $K$ is non-negative;
this is consistent with the condition of positivity of the estimated boundary function $f(\cdot)$.
The constraints (\ref{constr1S})--(\ref{constrPrime2S}) 
allow you to control the local properties of smoothness estimation
$\widehat{f}_N$ on the interval $[0, 1]$, which are used in the proofs.
It is interesting to note that the above described estimator
(\ref{estimator}),
(\ref{IPgoalS})--(\ref{constr2S}) may be treated as
the approximation to Maximum Likelihood Estimate related to the
estimation family~(\ref{estimator}); see {\sc Bouchard} {\it et al.}~\cite{BGIN}
for similar remarks.

\section{Main results}\label{secmain}

In the following theorem, the consistency and the convergence rate
of the estimator towards the true frontier is established with
respect to the $L_1$ norm over interval $[0,1]$.

\begin{Theo}\label{Th1}
Let the above mentioned assumptions A1, A2, and B1--B3 hold true and
the estimator parameter $C_{\alpha}\geq 8 f_{\max}$.
Moreover, let $h\to 0$ as $N\to\infty$ such that
\begin{equation}
\label{eq:limsup} \rho_{-} < \liminf_{N\to\infty} \frac{\log N}{N
h^{1+\beta}},\quad 
\limsup_{N\to\infty} \frac{\log N}{N h^{1+\beta}}
\leq \rho^+ < +\infty,
\end{equation}
where
\begin{equation}
\label{hypotheserho} \rho_- > \frac{f_{\max}}{L_{\beta}} \frac{C_X
K'_{\max}}{10 \times 3^\beta  K'_{\max}+ 3 C_{\beta}(K')}
\mbox{ and } C_X > 4 \frac{C_f}{f_{\min}}.
\end{equation}
Then, the LP-estimator
(\ref{estimator}), (\ref{IPgoalS})--(\ref{constr2S})
with kernel \eqref{estimatorsuitePlus3} for $(x,t)\in[-h,1+h]^2$
has the following
a.s.-properties:
\begin{equation}\label{L1rate}
\|\widehat{f}_N -f\|_{1} \leq \left( C_{12}(\beta)
[\rho_-]^{-\frac{\beta}{1+\beta}} + 2C_4(\beta)
[\rho^+]^{\frac{2}{1+\beta}}\right) \left(\frac{\log N}{N}\right)^\frac{\beta}{1+\beta}
 (1+o(1))
\end{equation}
asymptotically as $N\to\infty$ with constants
\begin{equation}\label{C12}
C_{12}(\beta) \triangleq 5 f_{\max}C_X\rho^+ K_{\max}
 + 10 \times 3^\beta L_{\beta} K_{\max}
+ 3 L_{\beta} C_{\beta}(K)
\end{equation}
and
\begin{equation}\label{ConstC4}
C_4(\beta) \triangleq 2 L_{\beta} \left(
\frac{2C_f}{f_{\min}L^{}_{\beta}\,}\right)^\frac{\beta}{1+\beta}
\left(\frac{1}{\rho_-} \right)^{\frac{2}{1+\beta}} +
 7 f_{\max} C_X K'_{\max} \left(\frac{2C_f}{f_{\min}L^{}_{\beta}\,}\right)^\frac{1}{1+\beta}\,.
\end{equation}
\end{Theo}
\medskip
Let us highlight that~(\ref{L1rate}) shows that $\widehat f_N$
reaches (up to a logarithmic factor) the minimax $L_1$ rate for
frontiers $f$ with H\"older continuous derivative, see {\sc Korostelev} \& {\sc
Tsybakov}~\cite{KorTsy2}, Theorem~4.1.1.

\section{Conclusions}

\textrm{The results obtained above straightly extend the approach, developed
in \cite{BGIN}, \cite{BGIN2005}, and \cite{GIN} under condition
$0<\beta\leq 1$, onto more smooth (and periodic) boundary functions
when the first derivative function is H\"older continuous with
exponent
 $\beta-1$,
and $1<\beta\leq 2$.
The estimation method itself for a boundary function,
like there in \cite{BGIN}, \cite{GIN},
reduces to a linear combination of sufficiently smooth kernel functions,
being centered at the sample points,
while weighting coefficients are defined by solving a linear programming (LP) problem
having minimized a sum of weighting coefficients under related constraints.
Note that the related LP problem changes in accordance with the value of  $\beta$
in a sense that the constraints composition in LP problem depends on degree of smoothness of boundary function:
when  $1<\beta\leq 2$ additional inequalities (\ref{constr1S}) under $X_i\neq X_j$ include into constraints
and the upper bound of second derivative of LP-estimate
 (\ref{constrPrime2S}) for all the sample points $X_i$.
Remind that under $0<\beta\leq 1$ the LP problem of \cite{GIN} contains the constraints
 (\ref{constr1S}) of this work only for $i=j$,
 i.e., inequalities of type $\widehat{f}_N(X_i) \geq Y_i$;
in addition, the upper bound on first derivative of LP-estimate (\ref{constrPrime2S})
is used for all the sample points $X_i$.
In particular it is evident that the transition of the parameter $\beta$ from interval $(0,1]$ to interval $(1,2]$
the number of constraints  of LP problem defining LP-estimate
increases stepwise which may be considered as
a certain fee for providing an optimal estimate under a smoother boundary function.
Note as well that all the values of  $\beta\in(0,2]$ give the error of estimation in $L_1$-norm
of the type $O\left((\log{N}/N)^{\beta/(1+\beta)}\right)$ a.s.
under the bandwidth selection of kernel function of type
$h\sim \left(\log{N}/N\right)^{1/(1+\beta)}$.
Finally, the authors hope to present an adaptive analog of the proposed method
at the given paper for the case of unknown a priori parameter smooth $\beta$.
\\
As a conclusion, the authors express their sincere thanks to the anonymous referees
for their critical comments contained the stimulating remarks and propositions.
}



\section{Appendix}


The proof of Theorem \ref{Th1} which is given in Subsection
\ref{PrTh1} is based on both upper and lower bounds derived in
Subsection~\ref{subsec:UB} and Subsection~\ref{subsec:LB}
respectively. When proving these bounds, we assume that the sequence
of the sample $X$-points $(X_i)_{i=\overline{1,N}}$ is already
increase ordered, without changing notation from $X_i$ to $X_{(i)}$
for the sake of simplicity, that is
\begin{equation}\label{Xorder}
X_i\leq X_{i+1}\,,\quad
\forall\, i \,.
\end{equation}
We essentially apply the uniform asymptotic bound $O(\log{N}/N)$ on
$\Delta X_i \triangleq X_i - X_{i-1}$ proved in auxiliary Lemma
\ref{CXlemma} (see also Lemma~A.2 in \cite{GIN}). Before that, we prove in
Subsection~\ref{subsec:proofpreli} some preliminary results.

\subsection{Proof of preliminary results}
\label{subsec:proofpreli}

\noindent\textbf{Proof of Lemma \ref{lem:IntegrKer}.} Under
assumptions of the Lemma, one may easily demonstrate that the kernel
definition (\ref{estimatorsuite}) ensures the equality
(\ref{eq:LemKerApprox1}).
Indeed, the LHS (\ref{eq:LemKerApprox1}) may be written from
(\ref{estimatorsuite}) for any $\,x\in[0,1]$ as follows:
\begin{eqnarray}
    \int_0^1 f(u) K_h(x,u)\,du &=&
    \label{eq:LemKerApprox1LHS1}
     \int_0^1 f(u) \frac{1}{h} K\left(\frac{x-u}{h}\right)\,du \\
     \label{eq:LemKerApprox1LHS2}
     &+&  \int_0^h f(u) \frac{1}{h} K\left(\frac{x-u-1}{h}\right)\,du \\
     \label{eq:LemKerApprox1LHS3}
     &+&  \int_{1-h}^1 f(u) \frac{1}{h} K\left(\frac{x-u+1}{h}\right)\,du
\end{eqnarray}
and, by changing the variables in the integrals
(\ref{eq:LemKerApprox1LHS1})--(\ref{eq:LemKerApprox1LHS3}), using
the $1$-periodicity assumption and the compact support assumption
B1, one may obtain
\begin{eqnarray}
    \int_0^1 f(u) K_h(x,u)\,du &=&
    \label{eq:LemKerApprox1LHS1a}
     \int_{(x-1)/h}^{x/h} f(x-hv) K(v)\,dv \\
     \label{eq:LemKerApprox1LHS2a}
     &+&  \int_{(x-1-h)/h}^{(x-1)/h}  f(x-hv) K(v)\,dv  \\
     \label{eq:LemKerApprox1LHS3a}
     &+&  \int_{x/h}^{(x+h)/h}  f(x-hv) K(v)\,dv \\
      \label{eq:LemKerApprox1LHS4a}
     &=& \int_{-1}^1  f(x-hv) K(v)\,dv\,.
\end{eqnarray}
Thus, (\ref{eq:LemKerApprox1}) is proved.
Now (\ref{eq:IntUKer}) follows directly from
 (\ref{eq:LemKerApprox1}) and assumption B2 for $f(x)\equiv 1$.
One obtains from similar arguments to
(\ref{eq:LemKerApprox1LHS1})--(\ref{eq:LemKerApprox1LHS4a}) that
\begin{eqnarray}\label{eq:LemKerApprox1LHSxu}
    \int_0^1 (x-u) K_h(x,u)\,du &=&
     h \int_{-1}^1  v\, K(v)\,dv \\
     \label{eq:LemKerApprox1LHS2axhv}
     &+&   \int_{(x-1-h)/h}^{(x-1)/h} K(v)\,dv  \\
     \label{eq:LemKerApprox1LHS3axhv}
     &-&  \int_{x/h}^{(x+h)/h}  K(v)\,dv\,.
\end{eqnarray}
Hence, condition $h\leq x\leq 1-h$  implies (\ref{eq:IntU-XdUKer}).
Finally, equalities (\ref{eq:IntXKer1}) and (\ref{eq:IntXUdxKer})
follow directly from the kernel definition (\ref{estimatorsuite})
and from the assumption B2:
     for all $\,u\in[0,1]$ one may easily verify that
\begin{equation}
     \int_0^1 K_h(x,u)\,dx = \int_{-1}^1 K(v)\,dv =
      1\,,
\end{equation}
and
\begin{equation}
     \int_0^1 (x-u) K_h(x,u)\,dx = h \int_{-1}^1 v\, K(v)\,dv =
      0\,.
\end{equation}
\CQFD\medskip

\noindent\textbf{Proof of Lemma \ref{LemSurfApprox}}
 ~is a straightforward consequence of
 (\ref{estimatorsuite}) and assumptions B1 and B2,
since
\begin{equation}\label{SurfApproxUBPr}
  \int_0^1 K_h(x,X_i)\, dx = 1\qquad\forall\,\, i=\overline{1,N}\,.
\end{equation}

\CQFD\medskip

\noindent\textbf{Proof of Lemma \ref{RegLemma1}.} For any
$x\in[0,1]$, one may write
\begin{eqnarray}\label{eq:MaxTermBoundPP1}
|\widehat{f}_N^{\;\prime}(x)| &\leq& \sum_{i=1}^N \alpha_i \left|
\frac{d}{dx}\,K_h( x,X_i) \right|
\\ \label{eq:MaxTermBound1PP1}
    &\leq& \sup_{u,v} \left|\frac{\partial}{\partial v}\,K_h(v,u) \right|
        \cdot     \sum_{i=1}^N \alpha_i \,\mathbf{1}\{|x-X_i|\leq h\}
\\ \label{eq:MaxTermBound2PP1}
    &\leq& 3 K'_{\max} C_\alpha h^{-1}\,,
\end{eqnarray}
where constraints (\ref{constrSumAlp}) are used and give
(\ref{eq:MaxTermBound2PP1}). This proves the Lemma. \CQFD\medskip

\noindent\textbf{Proof of Lemma \ref{RegLemma}.} We are to prove
(\ref{eq:RegLemmaP}). Remind that we assume (\ref{Xorder})
which imply $X_{-}=X_{0}$ and $X_{+}=X_{N+1}$ due to (\ref{eq:X0XN1}).
Consider the additional assumption
 (which holds true for all sufficiently large $N$)
that is
\begin{equation}\label{eq:hAddAssump}
    C_X^2\frac{\log N}{Nh}\leq 5\,
           \frac{L_{\beta} K''_{\max} }{ C_\alpha K''''_{\max}}           \,.
\end{equation}
By applying (\ref{constrPrime2S}) and auxiliary Lemma~\ref{CXlemma}
and Lemma~\ref{AuxLemma} we arrive at
\begin{eqnarray}\label{eq:ApplAuxLemmaPP}
&&    \max_{x\in[X_0,X_{N+1}]} |\widehat{f}_N^{\;\prime\prime}(x)|\\
&=& \max_{1\leq i\leq N+1}\, \max_{x\in[X_{i-1},X_i]}
|\widehat{f}_N^{\;\prime\prime}(x)|
\\ \label{eq:ApplAuxLemma1PP}
    &\leq&
 2\,L^{}_{\beta}\, K''_{\max}\, \frac{\log N}{Nh^3}
     +\frac{1}{8} \max_{1\leq i\leq N+1} \left[(X_i-X_{i-1})^2
             \max_{x\in[X_{i-1},X_i]} |\widehat{f}_N^{\;\prime\prime\prime\prime}(x)|
                                    \right]
\\ \label{eq:ApplAuxLemma2PP}
    &\leq&
 2\,L^{}_{\beta}\, K''_{\max}\, \frac{\log N}{Nh^3}
     +\frac{1}{8} \left(C_X \frac{\log{N}}{N}\right)^2
            \max_{x\in[X_0,X_{N+1}]} |\widehat{f}_N^{\;\prime\prime\prime\prime}(x)|\,,
\end{eqnarray}
with $C_X>4 C_f/f_{\min}$. The maximum term in
(\ref{eq:ApplAuxLemma2PP}) is bounded as follows: for any
$x\in[X_0,X_{N+1}]$\,,
\begin{eqnarray}\label{eq:MaxTermBoundPP}
|\widehat{f}_N^{\;\prime\prime\prime\prime}(x)|
&\leq& \sum_{i=1}^N \alpha_i 
\left| \frac{d^4}{dx^4}\,K_h( x,X_i) \right|
\\ \label{eq:MaxTermBound1PP}
    &\leq& \sup_{u,v} \left|\frac{\partial^4}{\partial v^4}\,K_h(v,u) \right|
        \cdot     \sum_{i=1}^N \alpha_i \,\mathbf{1}\{|x-X_i|\leq h\}
\\ \label{eq:MaxTermBound2PP}
    &\leq& 3 K''''_{\max} C_\alpha h^{-4}\,,
\end{eqnarray}
since
\begin{equation}\label{eq:LK2primHbPP}
 \sup_{u,v} \left|\frac{\partial^4}{\partial v^4}\,K_h(v,u) \right|
 \leq K''''_{\max} h^{-5}.
\end{equation}
Substituting (\ref{eq:MaxTermBoundPP}) and (\ref{eq:MaxTermBound2PP})
into (\ref{eq:ApplAuxLemma2PP}) and using (\ref{eq:hAddAssump})
yield
\begin{eqnarray}
    \max_{x\in[X_0,X_{N+1}]} |\widehat{f}_N^{\;\prime\prime}(x)|
    &\leq&
 2\,L^{}_{\beta}\, K''_{\max}\, \frac{\log N}{Nh^3}
     +  \frac{3}{8}\,K''''_{\max}  C_{\alpha}
     \left(C_X \frac{\log{N}}{N h^3}\right)^2   h^2\\
& \leq &
 4\,L^{}_{\beta}\, K''_{\max}\, \frac{\log N}{Nh^3}\,.
\end{eqnarray}
The result follows. \CQFD\medskip


\subsection{Upper bound for $\widehat{f}_N$ in terms of $J^*_P$}
\label{subsec:UB}

\begin{lemma}\label{Lm1Th1}
Let the assumptions of Theorem \ref{Th1} hold true.
Then for any 
\begin{equation}\label{GammaCondi}
    \gamma > \left(1+\frac{1}{\beta}\right) L^{}_{\beta}
    C_{\beta}(K)
     + f_{\max}\, (5K_{\max}+K'_{\max}) C_X \rho^+
\end{equation}
where $C_X> 4 C_f/f_{\min}$ and parameter $\rho^+$, meeting (\ref{eq:limsup}),
and for almost all $\omega\in\Omega$ there exist finite numbers
$N_1=N_1(\omega,\gamma)$ such that for all $N\geq N_1$ the LP
 (\ref{IPgoalS})--(\ref{constr2S})
 is solvable and
\begin{equation}\label{UB}
  J^*_P \leq C_f +\gamma h^\beta\,.
\end{equation}
\end{lemma}
\medskip

\noindent\textbf{Proof of Lemma \ref{Lm1Th1}.}
Recall that $\Delta X_{i}= X_i -X_{i-1}>0$ a.s. for all $i$
due to condition (\ref{Xorder}).
Consider arbitrary $N\geq
N_0(\omega)$ with $N_0(\omega)$ from Lemma \ref{CXlemma}. Introduce
function $f_\gamma(u)=f(u)+\gamma h^\beta$ with parameter $\gamma>0$
and pseudo-estimates
\begin{equation}\label{pseudoestimates}
\widetilde{\alpha}_i
 = \sum_{k=-1}^{1} a_{i,k} \int_{X_{i+k-1}}^{X_{i+k}} f_\gamma(u)\,du
 \,,\quad i=\overline{1,N}\,,
\end{equation}
where
\begin{equation}\label{aik}
  a_{i,k} = f_\gamma(X_{i+k}) \left(\int_{X_{i+k-1}}^{X_{i+k}} f_\gamma(u)\,du\right)^{-1}
                         \int_{X_{i+k-1}}^{X_{i+k}} b_{i,-k}(u)\,du\,,
\end{equation}
and functions $b_{i,k}(\cdot)$, $k=-1,0,1$, represent the coefficients
of the 2nd order Lagrange interpolation polynomial for the interval
defined by three successive points $X_{i-k-1}<X_{i-k}<X_{i-k+1}$,
i.e.,
\begin{equation}\label{bik}
  b_{i,k}(u) =\frac{\prod_{j=-1,\,j\neq k}^{j=1}(u-X_{i-k+j})}{\prod_{j=-1,\,j\neq k}^{j=1}(X_{i}-X_{i-k+j})}.
\end{equation}
For any 3 times continuously differentiable
function $g:[0,1]\to\mathbb{R}$ and for all $u\in[X_{i-1},X_{i+1}]$, the interpolation error is bounded by
\begin{eqnarray}\label{IntError1}
 \left| \sum_{k=-1}^{1} b_{i+k,k}(u) g(X_{i+k}) - g(u)
 \right|
 &\leq& 
 \max_{u\in [X_{i-1},X_{i+1}]} \left|\frac{g{'''}(u)}{6}\! \prod_{j=-1}^{1}\![u-X_{i+j}]\right|
  \\ \label{IntError3}
  &\leq& \frac{1}{9\sqrt{3}}\!
  \left[\max_{1\leq{i}\leq{N+1}}\!\Delta X_{i}\right]^3\!\! \max_{u\in[0,1]} \left|g{'''}(u)\right|\!.
\end{eqnarray}

\noindent $\mathrm{1}$. First, we prove constraints (\ref{constr1S})
under $\alpha_i = \widetilde{\alpha}_i$, $i=\overline{1,N}$. For
arbitrary $x\in[0,1]$,
\begin{eqnarray}\label{ProveConstraints1}
~~~~~~~~~~~~\widetilde{f}_N(x) &\triangleq&
     \displaystyle \sum_{i=1}^N \widetilde{\alpha}_i K_h(x,X_i)
     \,=\,\sum_{k=-1}^{1} \sum_{i=1}^N a_{i,k} \int_{X_{i+k-1}}^{X_{i+k}} f_\gamma(u)\,du\, K_h(x,X_i)
 \\ \label{ProveConstraints11}
&=& \displaystyle \sum_{k=-1}^{1} \sum_{i=1+k}^{N+k} a_{i-k,k}
\int_{X_{i-1}}^{X_{i}} f_\gamma(u)\,du\, K_h(x,X_{i-k})
 \\ \label{ProveConstraints11a}
&=& \displaystyle \sum_{i=1}^{N} \int_{X_{i-1}}^{X_{i}}
        f_\gamma(u)\,du\,\left[ a_{i,0}K_h(x,X_i)
                                + a_{i+1,-1}K_h(x,X_{i+1})
                                \right.
    \\ \label{ProveConstraints11b}
&& \phantom{\;=\; \displaystyle \sum_{i=1}^{N}
\int\limits_{X_{i-1}}^{X_i}
        f_\gamma(u)\,du\,\left(\right.}
                                \left.
                                + a_{i-1,1}K_h(x,X_{i-1})
                            \right]
 \\ \label{ProveConstraints12}
    && +\, a_{1,-1}K_h(x,X_1) \int\limits_{X_{-1}}^{X_0} f_\gamma(u)\,du -  a_{0,1}K_h(x,X_0) \int\limits_{X_0}^{X_1} f_\gamma(u)\,du
 \\ \label{ProveConstraints13}
    && +\, a_{N,1}K_h(x,X_N) \int\limits_{X_N}^{X_{N+1}} f_\gamma(u)\,du - a_{N+1,-1}K_h(x,X_{N+1}) \int\limits_{X_{N-1}}^{X_{N}} f_\gamma(u)\,du
 \\ \label{ProveConstraints14}
&=& \int\limits_{X_0}^{X_N} f_\gamma(u)\, K_h(x,u)\,du
 \\ \label{ProveConstraints15}
    && +\displaystyle \sum_{i=1}^{N} \int\limits_{X_{i-1}}^{X_{i}}\!\!
        f_\gamma(u) \left( \sum_{k=-1}^{1} a_{i+k,-k} K_h(x,X_{i+k}) -K_h(x,u)\!\!
                    \right)\!du\,,
\end{eqnarray}
since the sum of the terms
(\ref{ProveConstraints12})--(\ref{ProveConstraints13}) equals zero
due to the $1$-periodicity of function $f$ and the definition
(\ref{estimatorsuite})
 of the kernel $K_h$\,.
 Particularly, one can verify that $a_{1,-1}=a_{N+1,-1}$\,, $a_{0,1}=a_{N,1}$\,,
 $K_h(x,X_1)\equiv K_h(x,X_{N+1})$\,, and $K_h(x,X_0)\equiv K_h(x,X_N)$\,.
Now we separately bound each of the summands
 (\ref{ProveConstraints14})--(\ref{ProveConstraints15}) from below.

 Due to the kernel definition (\ref{estimatorsuite})
 and the conditions A1, B1, and B3,
 the main term (\ref{ProveConstraints14})  is transformed and bounded as follows:
\begin{eqnarray}
\label{ProveConstraints14sep}
 ~~~~~~~~\int\limits_{X_0}^{X_N} f_\gamma(u)\, K_h(x,u)\,du
  &=& \int\limits_{0}^{1} f_\gamma(u)\, K_h(x,u)\,du
 \\
 \label{ProveConstraints14sep0aa}
&&   + \int\limits_{X_0}^{0} f_\gamma(u)\,K_h(x,u)\,du
   - \int\limits_{X_N}^{1} f_\gamma(u)\,K_h(x,u)\,du
 \\
 \label{ProveConstraints14sep0}
  &=& \int\limits_{0}^{1} f(u)\, K_h(x,u)\,du
  + \gamma h^\beta
\end{eqnarray}
since the difference of two integrals in
(\ref{ProveConstraints14sep0aa}) vanishes due to periodicity assumption A1
and due to (\ref{eq:X0XN1}) 
and 
 (\ref{estimatorsuitePlus3}).
We continue
the integral in (\ref{ProveConstraints14sep0}) by Lemma~\ref{lem:IntegrKer}
as follows:
\begin{eqnarray}
\label{ProveConstraints14sepF}
 ~~~~~~~\int\limits_{0}^{1}\! f_\gamma(u)\, K_h(x,u)\,du
  &=& \int\limits_{-1}^{1}\! f(x-h{}t)\,K(t)\,dt
 +\, \gamma h^\beta
 \\
 \label{ProveConstraints14sepFb}
   &=& f(x)  +\, \gamma h^\beta
\\
 \label{ProveConstraints14sepF1b}
   &&
    +\int\limits_{-1}^{1} [f(x-h{}t) -f(x) -f'(x)(-h{}t)]\,K(t)\,dt
 \\
 \label{ProveConstraints14sepF1a}
  &\geq&
 f(x)
 + \left[\gamma  - \frac{L_{\beta}}{\beta}\, C_{\beta}(K) \right]\! h^{\beta}
 \,.
\end{eqnarray}
Notice that Lemma~\ref{CXlemma} as well as the
definitions 
for $C_{\beta}(\cdot)$ and
 $C_X$ have been used in
(\ref{ProveConstraints14sepF1b})--(\ref{ProveConstraints14sepF1a}).

The $i$-th summand from (\ref{ProveConstraints15}) which is denoted
below by $(\ref{ProveConstraints15})_i$ is decomposed and then
bounded basing on the 2nd order Lagrange interpolation with the
error upper bound (\ref{IntError1})--(\ref{IntError3}) being applied
for $g(u)=K_h(x,u)$ as follows:
\begin{eqnarray}\label{ProveConstraints15sep1st}
~~~~~~~(\ref{ProveConstraints15})_i
&\triangleq&\int\limits_{X_{i-1}}^{X_{i}}\!\!
        f_\gamma(u) \left( \sum_{k=-1}^{1} a_{i+k,-k} K_h(x,X_{i+k}) -K_h(x,u)\!\!
                    \right)\!du
\\ \nonumber && [\,\textrm{by~applying~definition~(\ref{aik})}\,]
 \\ \label{ProveConstraints15sep11st}
&=& f_\gamma(X_i)
         \sum_{k=-1}^{1}\int\limits_{X_{i-1}}^{X_{i}} b_{i+k,k}(u) K_h(x,X_{i+k})\,du
          -\int\limits_{X_{i-1}}^{X_{i}} f_\gamma(u) K_h(x,u)\,du
 \\ \label{ProveConstraints15sep21st}
&=& f_\gamma(X_i)
         \int\limits_{X_{i-1}}^{X_{i}} \left(\sum_{k=-1}^{1} b_{i+k,k}(u) K_h(x,X_{i+k})
                                              - K_h(x,u) \right) du
 \\ \label{ProveConstraints15sep21st1}
&&          -\int\limits_{X_{i-1}}^{X_{i}} \left(f_\gamma(u)
-f_\gamma(X_i) \right) K_h(x,u)\,du \,.
\end{eqnarray}
So, we apply Lemma~\ref{CXlemma} and the upper bound on the interpolation in
(\ref{IntError1})--(\ref{IntError3}):
\begin{eqnarray}\label{ProveConstraints15sep}
~~~~~(\ref{ProveConstraints15})_i &\geq& -2f_{\max} \, \frac{(\max\Delta
X_i)^3}{9\sqrt{3}} \max_{u\in[0,1]}
 \left|\frac{\partial^3 K_h(x,u)}{\partial u^3}\right|
 \,\Delta X_i\,
 \mathbf{1}\{|x-X_i|\leq 2h\}
 \\ \label{ProveConstraints15sep4}
 && - L^{}_{f}\,(\Delta X_i) \int_{X_{i-1}}^{X_i} K_h(x,u)\,du
 \\ \label{ProveConstraints15sep5}
&\geq& -\left( C_X \frac{\log{N}}{N}\right)^3
 \frac{2f_{\max}\,L_{K^{\prime\prime}}}{9\sqrt{3}h^4} \,\Delta X_i\,\mathbf{1}\{|x-X_i|\leq 2h\}
 \\ \label{ProveConstraints15sep5a}
&&  - L^{}_{f}\, C_X \frac{\log{N}}{N} \int_{X_{i-1}}^{X_i}
K_h(x,u)\,du.
\end{eqnarray}
Moreover, from Lemma~\ref{CXlemma}, it follows that
$$
\sum_{i=1}^{N} \Delta X_i \mathbf{1}\{|x-X_i|\leq 2h\} \leq 4h +
\frac{C_X \log N}{N}\,.
$$
Summing up by $i=\overline{1,N}$ we arrive at the bound for
the sum (\ref{ProveConstraints15}) as follows:
\begin{eqnarray}\label{ProveConstraints15sum}
[(\ref{ProveConstraints15})] &=& \sum_{i=1}^{N}
[(\ref{ProveConstraints15})_i]
\\ \label{ProveConstraints15sum51}
&\geq& -\left[ C_X \frac{\log{N}}{N}\right]^3\!
 \frac{f_{\max} L_{K^{\prime\prime}}}{4.5\sqrt{3}h^4} \left[4h + \frac{C_X \log N}{N}\right]
 \label{ProveConstraints15sum52}
 \! - L^{}_{f}  C_X \frac{\log{N}}{N} .
\end{eqnarray}
Thus, from (\ref{ProveConstraints1}),
(\ref{ProveConstraints14sep})--(\ref{ProveConstraints14sepF1a}),
(\ref{ProveConstraints15sep1st})--(\ref{ProveConstraints15sum51}),
and
(\ref{ProveConstraints14})--(\ref{ProveConstraints15}) it
follows for each $j=\overline{1,N}$ that
\begin{equation}\label{LastTerms}
\widetilde{f}_N(X_j) \geq f(X_j) + \delta_{0,N}
\end{equation}
with
\begin{eqnarray} \label{LastTerms1}
 \delta_{0,N} &\triangleq&
\left(\gamma - \frac{L^{}_{\beta}}{\beta}\,
    C_{\beta}(K)\right) h^\beta
\\ \label{LastTerms2}
&&
 -\left( C_X \frac{\log{N}}{Nh}\right)^3
        \frac{f_{\max}\,L_{K^{\,\prime\prime}}}{\sqrt{3}}
 - L^{}_{f}\, C_X \frac{\log{N}}{N}
 > 0
\end{eqnarray}
for sufficiently large $N\geq N_0(\omega)$ when
the following inequality holds true:
\begin{equation}\label{eq:LastCond}
 \gamma - \frac{L^{}_{\beta}}{\beta}\, C_{\beta}(K)
 \geq
 C_X\frac{\log N}{Nh^{1+\beta}}
  \left(
    \left( C_X \frac{\log{N}}{Nh}\right)^2
        \frac{f_{\max}\,L_{K^{\,\prime\prime}}}{\sqrt{3}}
    + L^{}_{f} h
  \right)\,.
\end{equation}
\medskip

\noindent $\mathrm{1^{\,'}}$.
%
Similarly, for arbitrary $x\in[0,1]$,
we now have to 
estimate the derivative value
\begin{equation}\label{ProveConstrPrime2S}
\widetilde{f}_N^{\;\prime}(x) =
     \displaystyle \sum_{i=1}^N \widetilde{\alpha}_i \,\frac{d}{dx}\,K_h(x,X_i) =
     \displaystyle \sum_{i=1}^N \widetilde{\alpha}_i \,\widetilde{K}_h(x,X_i),
\end{equation}
similarly to the arguments
(\ref{ProveConstraints1})--(\ref{ProveConstraints15sum52}). Here
\begin{equation}\label{DerivKerC2S}
\widetilde{K}_h(x,u) \triangleq \frac{\partial}{\partial
x}\,K_h(x,u)
 \end{equation}
with the following upper bound (see (\ref{estimatorsuite})\,)
\begin{equation}\label{TildeKerUB2S}
\left| \widetilde{K}_h(x,u) \right| \leq h^{-2} \max_{x}
\left|K'\left(\frac{x-u}{h}\right)\right| =  h^{-2} K'_{\max} \,.
\end{equation}
Hence,
one may repeat the arguments of (\ref{ProveConstraints11})--(\ref{ProveConstraints15}) 
by changing $K_h$ for $\widetilde{K}_h$, and, in particular,
equations (\ref{ProveConstraints1}),
(\ref{ProveConstraints14})--(\ref{ProveConstraints15}) give
\begin{eqnarray}\label{Der1Constraints1}
~~~~~\widetilde{f}_N^{\;\prime}(x) &= &\int\limits_{X_0}^{X_N}
f_\gamma(u)\, \widetilde{K}_h(x,u)\,du
 \\ \label{Der1Constraints15}
   &&
 +\displaystyle \sum_{i=1}^{N} \int\limits_{X_{i-1}}^{X_{i}}\!\!
        f_\gamma(u) \left( \sum_{k=-1}^{1} a_{i+k,-k} \widetilde{K}_h(x,X_{i+k}) - \widetilde{K}_h(x,u)\!\!
                    \right)\!du\,.
\end{eqnarray}
Therefore, all the rates from
(\ref{ProveConstraints15sep})--(\ref{LastTerms2})
should be divided by $h$, while the 
value of the main term of decomposition, due to conditions A1,
B1--B3,  as well as kernel representation
(\ref{estimatorsuitePlus3})
 is 
expressed as follows:
\begin{eqnarray}\label{ProveConstrPrimeSmain2S}
 \int\limits_{X_0}^{X_N} f_\gamma(u)\,\widetilde{K}_h(x,u) \,du
&=& \frac{1}{h}\int\limits_{-1}^{1} f(x-h{}t)\,K'(t)\,dt
 \\
 \label{ProCon14sepFbS}
&=& f'(x) +  \int\limits_{-1}^{1} [f'(x-h{}t) -f'(x) ]\,K(t)\,dt
 \\
 \label{ProCon14sepFbSs}
&=& f'(x)+ \delta_{1,N}
\end{eqnarray}
with
\begin{eqnarray}
\label{eq:delta4primS} |\delta_{1,N}| &\leq&  L_{\beta}
C_{\beta}(K) h^{\beta-1}\,;
\end{eqnarray}
cf. (\ref{ProveConstraints14sep})--(\ref{ProveConstraints14sepF1a}).
Furthermore,  the summation of integrals in
(\ref{Der1Constraints15}) gives the bound
\begin{equation}\label{Der1ConstrPrimeSmain}
    |(\ref{Der1Constraints15})| \leq \left( C_X
\frac{\log{N}}{N}\right)^3
 \frac{2f_{\max}\,L_{K^{\prime\prime\prime}}}{9\sqrt{3}h^4} \left(4 + \frac{C_X \log N}{Nh}\right)
+ L^{}_{f}\,  C_X \frac{\log{N}}{Nh}
 \end{equation}
instead of
(\ref{ProveConstraints15sep1st})--(\ref{ProveConstraints15sum52}).
Thus, by taking
(\ref{eq:delta4primS}) into
account, for sufficiently large $N\geq N_0(\omega)$ and for each
$x\in[0,1]$ we arrive at
\begin{eqnarray}\label{ProveConstrPrimeS12S}
\left|\widetilde{f}_N^{\;\prime}(x) - f'(x)\right| &\leq&
|\delta_{1,N}|
+\, O\left(\frac{\log^3 N}{N^3h^{4}}\right)
 +\, O\left( \frac{\log N}{Nh} \right) 
 \,.
\end{eqnarray}
\medskip

\noindent $\mathrm{1^{\,''}}$. So, we now take
(\ref{LastTerms})--(\ref{LastTerms1}) and
(\ref{ProveConstrPrimeS12S}) into account in order to prove
constraints (\ref{constr1S}): for any $|X_j - X_i|\leq h$, this yields
\begin{eqnarray*}
  \widehat{f}_N(X_i)
                    +\left(X_j - X_i\right)\widehat{f}^{\;\prime}_N(X_i)
 &\geq& f(X_i) + \left(X_j - X_i\right)f'(X_i) \\
   &&   +\, \delta_{0,N}
   - h|\delta_{2,N}|\\
   && +\, O\left(\frac{\log^3 N}{N^3h^{3}}\right)
 + O\left( \frac{\log N}{N} \right)
 \\
  &\geq& Y_j + \delta_{3,N}
\end{eqnarray*}
where (recalling that $\delta_{0,N}$ is defined in
(\ref{LastTerms1}))
\begin{eqnarray*}
    \delta_{3,N} &\triangleq& \delta_{0,N}
    - L_{\beta} C_{\beta}(K) h^{\beta}
    +\, O\left(\frac{\log^3 N}{N^3h^{3}}\right)
 +\, O\left( \frac{\log N}{N} \right)
 \\
 &=& \left(\gamma - \left(1+\frac{1}{\beta}\right) L^{}_{\beta}
    C_{\beta}(K)\right) h^\beta
\\ 
&&
 +\, O\left( \frac{\log^3{N}}{N^3 h^3}\right)
 +\, O\left( \frac{\log N}{N} \right)
\end{eqnarray*}
being positive for sufficiently large $N\geq N_0(\omega)$ when both
inequalities
(\ref{eq:LastCond}) hold
true and, additionally,
\begin{equation}\label{eq:gamma}
    \gamma > \left(1+\frac{1}{\beta}\right) L^{}_{\beta}
    C_{\beta}(K)
     \,.
\end{equation}
Notice, that inequality (\ref{eq:gamma}) implies
(\ref{eq:LastCond}).
\medskip

\noindent $\mathrm{2^{\,''}}$.  Similarly, constraints
(\ref{constrPrime2S}) hold true under $\alpha_i =
\widetilde{\alpha}_i$, $i=\overline{1,N}$. Indeed, for arbitrary
$x\in[0,1]$, we now have to bound  the absolute value of
\begin{equation}\label{ProveConstrPrime2}
\widetilde{f}_N^{\;\prime\prime}(x) =
     \displaystyle \sum_{i=1}^N \widetilde{\alpha}_i \,\frac{d^2}{dx^2}\,K_h(x,X_i) =
     \displaystyle \sum_{i=1}^N \widetilde{\alpha}_i \,\widetilde{\widetilde{K}}_h(x,X_i)
\end{equation}
instead of (\ref{ProveConstraints1}). Here
\begin{equation}\label{DerivKerC2}
\widetilde{\widetilde{K}}_h(x,u) \triangleq
\frac{\partial^2}{\partial x^2}\,K_h(x,u)
 \end{equation}
with the following upper bound (see (\ref{estimatorsuite})\,)
\begin{equation}\label{TildeKerUB2}
\left| \widetilde{\widetilde{K}}_h(x,u) \right| \leq h^{-3} \max_{x}
\left|K''\left(\frac{x-u}{h}\right)\right| =  h^{-3} K''_{\max} \,.
\end{equation}
Hence,
one may repeat the arguments of (\ref{ProveConstraints11})--(\ref{ProveConstraints15}) 
by changing $K_h$ for $\widetilde{K}_h$. Therefore, all the rates
from (\ref{ProveConstraints15sep})--(\ref{LastTerms1}) should be
divided by $h^2$, while the absolute value of the main term of
decomposition, due to conditions B1--B3,
 is bounded as follows:
\begin{equation}\label{ProveConstrPrimeSmain2}
    \left|\int\limits_{X_0}^{X_N} f_\gamma(u)\,\widetilde{\widetilde{K}}_h(x,u) \,du \right|
= \frac{1}{h^2}\left|\int\limits_{-1}^{1} f(x-h{}t)\,K''(t)\,dt
\right|
\end{equation}
\begin{eqnarray}
 \label{ProCon14sepFb}
&\leq& \frac{1}{h^2}\left|\int\limits_{-1}^{1} [f(x-h{}t) -f(x)
-f'(x)(-h{}t)]\,K''(t)\,dt
    \right|
 \\ \label{ProveConstrPrimeSmain12}
&\leq& \frac{2\,L^{}_{\beta}}{\beta(\beta+1)} \,
K''_{\max}\,h^{\beta-2}
\,,
\end{eqnarray}
instead of (\ref{ProveConstraints14sep})--(\ref{ProveConstraints14sepF1a}). 
Thus, for sufficiently large $N\geq N_0(\omega)$ and for each $X_j$
we arrive at
\begin{eqnarray}\label{ProveConstrPrimeS12}
\left|\widetilde{f}_N^{\;\prime\prime}(X_j)\right| &\leq&
\frac{2\,L^{}_{\beta}\, K''_{\max}}{\beta(\beta+1)\,h^{2-\beta}}
 +\, O\left(
\frac{\log^3 N}{N^3h^{5}}\right)
 +\, O\left( \frac{\log N}{Nh^2} \right) 
\\ \label{ProveConstrPrimeS12b}
 &\leq&
\frac{3\,L^{}_{\beta}}{\beta\,\rho_{-}}\,
 K''_{\max}\, \frac{\log N}{Nh^3}.
\end{eqnarray}
Namely, inequality (\ref{ProveConstrPrimeS12b}) holds true almost
surely for all those $N\geq N_0(\omega)$ such that inequalities
(\ref{eq:AddAssumptions}) hold true
 and
\begin{eqnarray}\label{eq:AddAssumptions22}
\frac{L^{}_{\beta}}{\beta\,C_X}
 \left( \frac{1}{\rho_{-}} -\frac{h^{1+\beta}N}{\log N}
  \right)
&\geq&  \left( C_X \frac{\log{N}}{Nh}\right)^2
        \frac{f_{\max}\,L_{K^{\prime\prime\prime\prime}}}{\sqrt{3}\,K''_{\max}}
 + 2L^{}_{f}\, h
 \,.
\end{eqnarray}
\medskip

\noindent 3. Finally, the constraints (\ref{constrSumAlp}) with
\begin{equation}\label{CalphaValue}
C_\alpha\geq 
4 f_{\max}
\end{equation}
also hold true under $\alpha_i = \widetilde{\alpha}_i$,
$i=\overline{1,N}$. Indeed, by Lemma \ref{CXlemma} the following
inequalities hold \textrm{a.s.} for all $N\geq N_0(\omega)$ and for
each $j=\overline{1,m_h}$, where $m_h=\lfloor{h^{-1}}\rfloor$ :
\begin{eqnarray}\label{eq:FinalConstrCheck}
&& \sum_{i=1}^N \widetilde{\alpha}_i \, \mathbf{1}\{ (j-1)/m_h\leq X_i
< j/m_h\}
\\ \label{eq:FinalConstrCheck0}
&\leq& (f_{\max}+\gamma h^\beta) \left( 1/m_h + 2C_X
\frac{\log{N}}{N}\right)
\\
    &\leq& 
    4 f_{\max} h,
\end{eqnarray}
under additional assumptions
\begin{equation}\label{eq:AddAssumptions}
f_{\max}\geq \gamma h^\beta\,,
 \quad h\geq 2C_X^{}\log N/N\,.
\end{equation}
Thus,
constraints (\ref{constrSumAlp}) are fulfilled under
(\ref{CalphaValue}) almost surely, for any sufficiently large $N$.
\medskip

\noindent 4. Since all $\widetilde{\alpha}_i\geq 0$, constraints
(\ref{constr2S}) hold true. Hence, vector
$(\widetilde{\alpha}_1,\dots,\widetilde{\alpha}_N)^T$ is the
admissible point for the  LP
 (\ref{IPgoalS})--(\ref{constr2S}).
Now
inequality (\ref{UB}) follows from Lemma~\ref{LemSurfApprox}.
\CQFD
\medskip

\subsection{Lower bound for estimate $\widehat{f}_N$}
\label{subsec:LB}

\begin{lemma}\label{Lm2Th1}
Under the assumptions of Theorem \ref{Th1}, for almost all
$\omega\in\Omega$ there exist finite numbers $N_2(\omega)$
such that for any $x\in [0,1]$ and for all $N\geq N_2(\omega)$
\begin{equation}\label{LB}
  \widehat{f}_N(x) - f(x)
  \geq   - \frac{2L^{}_{\beta}}{\beta} h^{\beta}
            - \left( \frac{2 C_f}{f_{\min}} + 2\,L^{}_{\beta}\, K''_{\max} + 4C_{\alpha} K'_{\max} C_X \right)\,\frac{\log{N}}{Nh}
\end{equation}
assuming that
$
L_f h\leq C_{\alpha}K'_{\max}
$
and 
$C_X>4 C_f/f_{\min}$.
\end{lemma}
\medskip

%
\noindent\textbf{Proof of Lemma \ref{Lm2Th1}.} Let us take use of
Lemma~\ref{AuxPartiLemma} and its Corollary~\ref{CorAuxPartiLemma}
introducing
\begin{equation}\label{dely}
\delta_y \thicksim 
                \frac{\log{N}}{Nh}\,
                ,\quad \delta_x \thicksim h
\,.
\end{equation}
Thus, for any $N\geq N_6(\omega)$ and any $x\in[0,1]$ there exist
(with probability one)  integers $i_k\in\{1,\dots,N\}$,
$k=\overline{1,m_h}$, such that
\begin{equation}\label{epsuppose}
 |x-X_{i_k}| \leq \delta_x
 \end{equation}
and inequality (\ref{eq:LinAppWorks}) that is
\begin{equation}\label{deltay}
    \mathcal{L}_{x}f(X_{i_k})
    \leq Y_{i_{k}} + \delta_y +\,\frac{L_{\beta}}{2}\,\delta_x^{\beta}\,.
\end{equation}
Now, we put a point $x\in[0,1]$, find index $i_x\in\{1,\dots,N\}$
such that $|X_{i_x}-X_{i_k}|\leq h$ and
$$|x- X_{i_x}|\leq\max\{\Delta X_{i_{x}-1},\Delta X_{i_{x}}\}
\leq C_X\frac{\log N}{N};
$$
 then the estimation error at point $x$ can be decomposed as
\begin{eqnarray}\label{ferr01}
  f(x)-\widehat{f}_N(x) &=&  \left[f(x) - f(X_{i_x})\right]\\
 &+&
\label{ferr1}     \left[f(X_{i_x}) - \mathcal{L}_{X_{i_x}} f(X_{i_k})\right]\\
 &+&
\label{ferr2}        \left[\mathcal{L}_{X_{i_x}} f(X_{i_k}) - \mathcal{L}_{X_{i_x}} \widehat{f}_N(X_{i_k})\right]\\
 &+&
 \label{ferr3}        \left[\mathcal{L}_{X_{i_x}} \widehat{f}_N(X_{i_k}) - \widehat{f}_N(X_{i_k})\right]\\
 \label{ferr31}      &+&  \left[\widehat{f}_N(X_{i_k}) - \widehat{f}_N(x)\right] .
\end{eqnarray}
The first and the last decomposition components, i.e. RHS
(\ref{ferr01}) and (\ref{ferr31}),  can be similarly bounded by
using the proper Lipschitz constants as follows:
\begin{eqnarray}\label{eq:ferrComp01-31a}
    |f(x) - f(X_{i_x})|  &\leq& L_f |x - X_{i_x}|\,\leq\,L_f C_X\,\frac{\log N}{N}\,,
    \\
\label{eq:ferrComp01-31b}
    |\widehat{f}_N(X_{i_k}) - \widehat{f}_N(x)| &\leq& L_{\widehat{f}_N} |x - X_{i_x}| \,\leq\, L_{\widehat{f}_N}C_X \, \frac{\log N}{N}\,.
\end{eqnarray}
The similar decomposition components (\ref{ferr1}) and (\ref{ferr3})
may be bounded as follows:
\begin{eqnarray}\label{err1a}
  \left|f(X_{i_x})-\mathcal{L}_{X_{i_x}} f(X_{i_k})\right|
  &\leq&
  \frac{L^{}_{\beta}}{\beta}  \,  \left|X_{i_x}-X_{i_k} \right|^{\beta}
  \, \leq\,
  \frac{L^{}_{\beta}}{\beta}\, \delta_x^{\beta}\,,
\\
\label{err1b}
 \left| \mathcal{L}_{X_{i_x}} \widehat{f}_N(X_{i_k}) - \widehat{f}_N(X_{i_x}) \right|
  &\leq& \frac{L_{\widehat{f}^{\,\,\prime}_N}}{2}\,  \left|X_{i_x}-X_{i_k} \right|^2
   \,\leq\, \frac{L_{\widehat{f}^{\,\,\prime}_N}}{2}\, \delta^2_x,
\end{eqnarray}
with Lipschitz constant $L_{\widehat{f}^{\,\,\prime}_N}$ for the
derivative estimator function $\widehat{f}^{\,\,\prime}_N(x)$.
Finally, we bound the central decomposition component (\ref{ferr2})
by applying Corollary~\ref{CorAuxPartiLemma} of
Lemma~\ref{AuxPartiLemma} with
\begin{equation}\label{eq:delXdelY}
    \delta_x =h,\quad  \delta_y = \frac{2C_f}{f_{\min}}\, \frac{\log N}{Nh},
\end{equation}
and using the estimator constraints (\ref{constr1S}); we obtain
\begin{eqnarray}\label{eq:CentralErrA}
     \mathcal{L}_{X_{i_x}} f(X_{i_k}) - \mathcal{L}_{X_{i_x}} \widehat{f}_N(X_{i_k})
     &\leq&
         Y_{i_{k}} +  \delta_y +\,\frac{L_{\beta}}{2}\,\delta_x^{\beta} - Y_{i_{k}}
         \\ \label{eq:CentralErrB}
        &=& \frac{2C_f}{f_{\min}}\, \frac{\log N}{Nh} +\,\frac{L_{\beta}}{2}\,h^{\beta}\,.
\end{eqnarray}
Therefore, equations (\ref{ferr01})--(\ref{eq:CentralErrB}) and
Lemmas~\ref{RegLemma1} and \ref{RegLemma} lead to the lower bound
\begin{eqnarray}\label{LBestim}
&& \widehat{f}_N(x) -f(x)
\\ \label{LBestim0}
&\geq& - \left( \frac{2 C_f}{f_{\min}}\,
\frac{\log N}{Nh}
 + \frac{2L^{}_{\beta}}{\beta} h^{\beta}
 +\frac{L_{\widehat{f}^{\,\,\prime}_N}}{2} h^2
+(L_f + L_{\widehat{f}_N}) C_X \frac{\log N}{N}
     \right)
\\ \label{LBestim1}
&\geq &  - \frac{2L^{}_{\beta}}{\beta} h^{\beta}
            - \left( \frac{2 C_f}{f_{\min}} + 2\,L^{}_{\beta}\, K''_{\max} + 4C_{\alpha} K'_{\max} C_X \right)\,\frac{\log{N}}{Nh}
\end{eqnarray}
assuming additionally that
$$
L_f h\leq C_{\alpha}K'_{\max}\,.
$$
Thus, the obtained lower bound holds true for any sufficiently large
$N$ (starting from random a.s. finite integer, which does not depend
on $x$).
Lemma \ref{Lm2Th1} is proved. \CQFD

\begin{Remark}\label{rem:OptimLowerBound}
The optimal order of the lower bound, proved in Lemma~\ref{Lm2Th1},
is attained by
\begin{equation}\label{eq:Optimh}
    h=h_1\left(\frac{\log N}{N}\right)^{\frac{1}{1+\beta}}
\end{equation}
when two terms in (\ref{LBestim1}) are balanced, and the lower bound
in (\ref{LBestim})--(\ref{LBestim1}) becomes
\begin{equation}\label{eq:h1LB}
\widehat{f}_N(x) -f(x) \geq  - C_{LB}(h_1)\left(\frac{\log
N}{N}\right)^{\frac{\beta}{1+\beta}}
\end{equation}
where constant
\begin{equation}\label{eq:Ch1}
    C_{LB}(h_1) =   \frac{2L^{}_{\beta}}{\beta} h_1^{\beta}
            + \left( \frac{2 C_f}{f_{\min}} + 2\,L^{}_{\beta}\, K''_{\max} + 4C_{\alpha} K'_{\max} C_X \right)\,\frac{1}{h_1}
\end{equation}
may be optimized by $h_1>0$. It is interesting to observe that four
last components of the estimation error decomposition
(\ref{ferr01})--(\ref{ferr31}) become of the same order while the
first one RHS (\ref{ferr01}) be negligible, see
(\ref{eq:ferrComp01-31a})--(\ref{eq:CentralErrB}).
\end{Remark}

\subsection{Proof of Theorem \ref{Th1}}\label{PrTh1}
 \medskip \noindent  1. Since $|u|=u-2u\mathbf{1}\{u<0\}$,
 the $L_1$-norm of estimation error can be expanded as
\begin{eqnarray}\label{L1norm1}
  \|\widehat{f}_N -f\|_{1} &=& \int_0^1 \left[\widehat{f}_N(x) -f(x)\right]\,dx
\\ \label{L1norm2}
& +&2\int_0^1 \left[f(x)-\widehat{f}_N(x) \right] \mathbf{1}\!
\left\{\widehat{f}_N(x)<f(x)\right\}\,dx.
\end{eqnarray}
 \medskip

 \noindent  2.  Applying Lemmas \ref{LemSurfApprox} and \ref{Lm1Th1} to the right hand side (\ref{L1norm1}) yields
\begin{equation}\label{useUB}
  \limsup_{N\to\infty}\, h^{-\beta}
    \left(\int_0^1 \left[\widehat{f}_N(x) -f(x)\right]\,dx\right)
    \leq \gamma
\quad\mathrm{a.s.}
\end{equation}
where $\gamma>0$ is large enough.
\medskip

\noindent  3. In order to obtain a similar result for the term
(\ref{L1norm2}), note that Lemma \ref{Lm2Th1}
implies
$$
\zeta_N(x,\omega) \triangleq \varepsilon_{LB}^{-1}(N)
\left[f(x)-\widehat{f}_N(x) \right]
\leq \mathrm{const}
<\infty \quad \mathrm{a.s.}
$$
uniformly with respect to both $x\in[0,1]$ and $N\geq N_2(\omega)$,
with
\begin{equation}\label{epsLB}
  \varepsilon_{LB}(N) \triangleq
            \mathrm{const}\,\frac{\log N}{Nh}
\end{equation}
with finite $\mathrm{const}>0$. Hence, one may apply Fatou's lemma,
taking into account that $u\mathbf{1}\{u>0\}$ is a continuous,
monotone function:
\begin{eqnarray}
  &&  \limsup_{N\to\infty}\, \varepsilon_{LB}^{-1}(N)
  \int_0^1 \left[f(x)-\widehat{f}_N(x) \right] \mathbf{1}\! \left\{\widehat{f}_N(x)<f(x)\right\}\,dx
   \\
 &\leq&
 \int_0^1 \limsup_{N\to\infty}\, \zeta_N(x,\omega)\, \mathbf{1}\!
    \left\{\zeta_N(x,\omega)>0\right\}\,dx
   \\
 &\leq& \mathrm{const}
  <\infty \quad \mathrm{a.s.}
\end{eqnarray}

\noindent 4. Finally, we put bandwidth $h$ from the balancing
assumption
$$
h^{\beta} \thicksim \frac{\log N}{Nh} \,.
$$
Thus, the obtained relations together with (\ref{L1norm1}) and
(\ref{L1norm2}) imply (\ref{L1rate}). Theorem~\ref{Th1} is proved.
\CQFD

The following results are quoted here for the sake of completeness.
\begin{lemma}[(Lemma~A.2 in \cite{GIN})]\label{CXlemma}
Let function $f:[0,1]\to\R$ meets the assumption A1 and sequence
$(X_i)_{i=\overline{1,N}}$ be obtained from an independent sample
with p.d.f. $f(x)/C_f$ by increase ordering (\ref{Xorder}), where
$C_f$ is defined by (\ref{eq:Cf}). Denote $X_0=0$ and $X_{N+1}=1$.
Then for any finite constant $C_X > 4C_f/f_{\min}$ there exist
almost surely finite number $N_0=N_0(\omega)$ such that
\begin{equation}\label{LemmaX}
\max_{i=\overline{1,N+1}} \Delta{X_i}  \leq C_X \frac{\log{N}}{N}
\quad \forall \;N\geq N_0
\end{equation}
with probability 1. For instance, one may fix constant $C_X$ as
follows:
\begin{equation}\label{eq:CXvalue}
C_X= 5f_{\max}/f_{\min}\,.
\end{equation}
\end{lemma}
\medskip

\begin{lemma}[(Lemma~A.3 in \cite{GIN})]\label{AuxPartiLemma}
Let random sample $\{(X_i,Y_i)\, |\,\,i=\overline{1,N}\}$ be defined
as in Section~\ref{subsecLPprobl}. Let sequence
$\delta_x=\delta_x(N)$ be positive,  and for some $\varepsilon>0$
\begin{equation}\label{eq:liminfN1eps}
\liminf_{N\to\infty}N^{1-\varepsilon}\delta_x>0\,.
\end{equation}
Define
\begin{equation}\label{eq:mdeltaDef}
m_{\delta}^{} \triangleq \min\{\mathrm{integer}\; m\,:\; m\geq
\delta_x^{-1} \}
\end{equation}
and assume a positive sequence $\delta_y=\delta_y(N)<f_{\min}$
meeting for all sufficiently large $N$ the inequality
\begin{equation}\label{eq:deltayDef}
\delta_y^{}\geq \kappa\, m_{\delta}^{}\frac{\log N}{N}\,,\quad
\mathrm{with}\quad\kappa> \frac{(2-\varepsilon)C_f}{f_{\min}}\,.
\end{equation}
Then, under the assumptions of Lemma~\ref{CXlemma}, with probability
1, there exists finite number $N_6(\omega)$ such that for any $N\geq
N_6(\omega)$ there is such a subset of points
$\left\{\left(X_{i_k},Y_{i_k}\right),\,
k=\overline{1,m_{\delta}^{}}\right\}$\, in the sample
$\left\{\left(X_{i},Y_{i}\right),\right.$ $\left.
i=\overline{1,N}\right\}$, that the following inequalities hold:
\begin{equation}\label{eq:AuxPartiLemma}
(k-1)/m_{\delta}^{}\leq X_{i_k}< k/m_{\delta}^{}\,,\quad
f(X_{i_k})-\delta_y^{}\leq Y_{i_k}\leq f(X_{i_k})\,.
\end{equation}
\end{lemma}

\begin{Coro}\label{CorAuxPartiLemma}
Let $\delta_x$ and $\delta_y$ meet the conditions of
Lemma~\ref{AuxPartiLemma}. Then, with probability~1, for any $N\geq
N_6(\omega)$ and any $x\in[0,1]$ there exists integer
$i_k\in\{1,\dots,N\}$ such that $|x-X_{i_k}|\leq \delta_x$ and $
f(X_{i_k})-\delta_y^{}\leq Y_{i_k}\leq f(X_{i_k})$\,. Furthermore,
if constant Lipschitz $L_{\beta}<\infty$ with $1<\beta\leq2$, then
\begin{equation}\label{eq:LinAppWorks}
    \mathcal{L}_{x}f(X_{i_k})
    \leq Y_{i_{k}} + \delta_y +\,\frac{L_{\beta}}{2}\,\delta_x^{\beta}\,.
\end{equation}

\end{Coro}

\medskip

\begin{lemma}[(Lemma~A.4 in \cite{GIN})]\label{AuxLemma}
Let function $g:[0,\Delta]\to\R$ be twice continuous differentiable,
$\Delta>0$. Then
\begin{equation}\label{eq:AuxLemma}
    \max_{x\in[0,\Delta]} |g(x)| \leq \max\{ |g(0)|, |g(\Delta)| \}
     +\frac{\Delta^2}{8} \max_{x\in[0,\Delta]} |g^{\,\prime\prime}(x)|
     \,.
\end{equation}
\end{lemma}

\bigskip

{Nazin A.V., }{Trapeznikov Institute of Control Sciences RAS,
65, Profsoyuznaya str.,
Moscow, Russia, }{nazine@ipu.ru}

{Girard S., }{LJK, Inria Grenoble Rh\^one-Alpes, \'equipe-projet Mistis,
Inovall\'ee, 655, av. de l'Europe, Montbonnot, 38334 Saint-Ismier cedex,
 France, }{stephane.girard@inria.fr}


\end{document}